\documentclass [12pt]{amsart}

\usepackage{amssymb,amsxtra,amsfonts}
\usepackage{amsmath, accents} %
\usepackage{graphics}

\usepackage[colorlinks]{hyperref}

\usepackage{verbatim}  
\usepackage{bbding}   

\newenvironment{ppb}[1]
{\ \!\!\!\!\!\!\!\!\!\!\!\!\!\!\!\!\!\!\!\!\!\!\!\!\!\!\!\!\!\!\!\!\!\!\!\!\!\!\!\! {\bf PPB-----------------------------------------------------------------------------PPB}\newline \tiny {#1}
\  \newline\normalsize\phantom{f}\!\!\!\!\!\!\!\!\!\!\!\!\!\!\!\!\!\!\!\!\!\!\!\!\!\!\!\!\!\!\!\!\!\!\!\!\!\!\!\! {\bf 
PPB-----------------------------------------------------------------------------PPB}\newline}{}

\newcounter{parag}

\long\def\pb #1*/{}

\def\reE@DeclareMathSymbol#1#2#3#4{%
    \let#1=\undefined
    \DeclareMathSymbol{#1}{#2}{#3}{#4}}
\DeclareSymbolFont{symbolsC}{U}{txsyc}{m}{n}
\SetSymbolFont{symbolsC}{bold}{U}{txsyc}{bx}{n}
\DeclareFontSubstitution{U}{txsyc}{m}{n}
\reE@DeclareMathSymbol{\strictiff}{\mathrel}{symbolsC}{76}

\newcommand\beq{\begin{equation}}
\newcommand\eeq{\end{equation}}
\newcommand\bal{\begin{align*}}
\newcommand\eal{\end{align*}}   
\newcommand\bmx{\left(\begin{matrix}}
\newcommand\emx{\end{matrix}\right)}
\newcommand\bsmx{\left(\begin{smallmatrix}}
\newcommand\esmx{\end{smallmatrix}\right)}
\newcommand\bmxnp{\begin{matrix}}
\newcommand\emxnp{\end{matrix}}
\newcommand\bsmxnp{\begin{smallmatrix}}
\newcommand\esmxnp{\end{smallmatrix}}

\newcommand{\bSi}{{\bf \Si}}

\newcommand{\monesp}{\;\!\!}   
\newcommand{\onto}{\twoheadrightarrow}

\newcommand{\spq}{/\!\!/}
\newcommand{\sqp}{\setminus\!\!\!\monesp\setminus\,}

\providecommand{\spqa}[1]{\underset{#1}{/\!\!/}}

\newcommand{\st}{\ \bigl\vert\ }

\def\part#1{\frac{\partial\phantom{q}}{\partial#1}}

\newcommand {\flb}{\lbrack\!\lbrack}
\newcommand {\frb}{\rbrack\!\rbrack}
\newcommand {\flp}{(\!(}
\newcommand {\frp}{)\!)}

\newcommand{\fus}{\circledast}

\newcommand{\MDR}{\mathcal{M}_{\text{\rm DR}}}

\newcommand{\MB}{\mathcal{M}_{\text{\rm B}}}

\newcommand{\MDol}{\mathcal{M}_{\text{\rm Dol}}}

\newcommand{\HH}{\text{\rm H}}

\newcommand{\Lie}{{\mathop{\rm Lie}}}

\newcommand{\res}{{\mathop{\rm Res}}}

\newcommand{\Prod}{\prod}

\newcommand{\tr}{{\mathop{\rm Tr}}}

\DeclareMathOperator{\Hom}{Hom}         

\newcommand{\SL}{{\mathop{\rm SL}}}

\newcommand{\GL}{{\mathop{\rm GL}}}

\newcommand*{\SU}{{\mathop{\rm SU}}}
\newcommand{\U}{{\rm {U}}}	
\newcommand{\SO}{{\mathop{\rm SO}}}
\newcommand{\Sp}{{\mathop{\rm Sp}}}
\newcommand{\Spin}{\mathop{\rm Spin}}

\newcommand{\irr}{{\rm irr}}

\DeclareMathOperator{\Rep}{\rm Rep}
\DeclareMathOperator{\End}{End}

\newcommand{\hk}{{hyperk\"ahler }}   

\newcommand{\ba}{{\bf a}}
\newcommand{\bA}{{\bf A}}

\newcommand{\bB}{{\bf B}}

\newcommand{\bLa}{{\bf \Lambda}}

\newcommand{\bM}{{\bf M}}

\newcommand{\bQ}{{\bf Q}}

\newcommand{\bS}{{\bf S}}

\DeclareSymbolFont{bbold}{U}{bbold}{m}{n}
\DeclareSymbolFontAlphabet{\mathbbold}{bbold}

\newcommand{\IB}{\mathbb{B}}
\newcommand{\IC}{\mathbb{C}}
\newcommand{\ID}{\mathbb{D}}

\newcommand{\IH}{\mathbb{H}}

\newcommand{\IP}{\mathbb{P}}                                     
                           
\newcommand{\IR}{\mathbb{R}}                           
\newcommand{\IS}{\mathbb{S}}

\newcommand{\IZ}{\mathbb{Z}}

\newcommand{\cA}{\mathcal{A}}
\newcommand{\cB}{\mathcal{B}}
\newcommand{\Bun}{\text{Bun}}
\newcommand{\cC}{\mathcal{C}}

\newcommand{\cG}{\mathcal{G}}

\newcommand{\cL}{\mathcal{L}}

\newcommand{\cM}{\mathcal{M}}

\newcommand{\cO}{\mathcal{O}}

\newcommand{\cP}{\mathcal{P}}

\newcommand{\bcO}{\boldsymbol{\mathcal{O}}}

\newcommand{\bcM}{{\boldsymbol{\mathcal{M}}}}

\newcommand{\gM}{       \mathfrak{M}     }

\newcommand{\g}{       \mathfrak{g}     }

\newcommand{\lt}{\mathfrak{t}}
\newcommand{\lh}{\mathfrak{h}}

\renewcommand{\sl}{       \mathfrak{sl}     } 

\newcommand{\wt}{\widetilde}

\DeclareMathSymbol{\widehatsym}{\mathord}{largesymbols}{"62}
\newcommand\lowerwidehatsym{%
  \text{\smash{\raisebox{-1.3ex}{%
    $\widehatsym$}}}}
\newcommand\fixwidehat[1]{%
  \mathchoice
    {\accentset{\displaystyle\lowerwidehatsym}{#1}}
    {\accentset{\textstyle\lowerwidehatsym}{#1}}
    {\accentset{\scriptstyle\lowerwidehatsym}{#1}}
    {\accentset{\scriptscriptstyle\lowerwidehatsym}{#1}}
}

\newcommand{\wh}{\fixwidehat}

\newcommand{\Ga}{\Gamma}

\newcommand{\la}{\lambda}
\newcommand{\La}{\Lambda}

\newcommand{\si}{\sigma}

\newcommand{\Om}{\Omega}
\newcommand{\Si}{\Sigma}
\newcommand{\Th}{\Theta}
\renewcommand{\th}{\theta}

\makeatletter
 \newlength{\typesize}
\makeatother

\newlength{\vvoff}
\newlength{\hhoff}

\def\mapdown#1{\Big\downarrow
        \rlap{$\vcenter{\hbox{$\!\scriptstyle#1$}}$}}

\def\overcong#1{\smash{
        \mathop{\cong}\limits^{#1}}}

\def\underset#1#2{\ \smash{\mathop{ #2 }\limits_{#1}}\ }
\def\overcong#1{\smash{
        \mathop{\cong}\limits^{#1}}}

\newcommand{\pf}{\begin{bpf}}

\newcommand{\pfms}{\begin{bpfms}}
\newcommand{\epf}{\end{bpf}\hfill$\square$\\}           
\newcommand{\epfms}{\end{bpfms}\hfill$\square$\\}       

\newcommand{\idea}{\begin{bidea}}

\newcommand{\eidea}{\end{bidea}\hfill$\square$\\}           

\newcommand{\sk}{\begin{bsk}}    

\newcommand{\esk}{\end{bsk}\hfill$\square$\\}           
\newcommand{\sketch}{\begin{bsketch}}

\newcommand{\esketch}{\end{bsketch}\hfill$\square$\\}

\newtheorem {hypo}{\bf\hspace{-\parindent}Hypothesis}
\newtheorem {thm}[hypo]{Theorem}   

\newtheorem {cor}[hypo]{Corollary}
\newtheorem {lem}[hypo]{Lemma}

\newtheorem {defn}[hypo]{Definition}

\theoremstyle{remark}\newtheorem{rmk}[hypo]{Remark}

\begin{document}

\title[Wild character varieties, Hitchin systems and Dynkin diagrams]{Wild character varieties, meromorphic Hitchin systems and Dynkin diagrams}
\author{Philip Boalch}%

\begin{abstract}
The theory of Hitchin systems is something like a 
``global theory of Lie groups'', where one works 
over a Riemann surface rather than just at a point.
We'll describe how one can take this analogy a few steps further by attempting to make precise the class of rich geometric objects that appear in this story (including the non-compact case), and discuss their classification, outlining a theory of ``Dynkin diagrams'' as a step towards classifying some examples 
of such objects.
\end{abstract}

\maketitle

\section{The Lax Project}

We would like to try to classify integrable systems, upto isomorphism (or isogeny, deformation...).
For this we need a definition---here we will use the following:

\begin{defn}\label{def: acihs}
A finite dimensional complex algebraic integrable Hamiltonian system is a symplectic algebraic variety $M$ with a map $\chi:M\to \IH$ to an affine space $\IH$ (of half the dimension) such that generic fibres of $\chi$ are Lagrangian abelian varieties.
\end{defn}

This is close to the classical point of view of having $n$ independent Poisson commuting functions on a symplectic manifold of dimension $2n$.
One can consider a more general definition allowing generic fibres to be open parts of Lagrangian abelian varieties, but for the cases we look at there are natural compactifications of the fibres.
See \cite{vanh-aci} Ch. 5 for a discussion of some other possible definitions.

Of course this is a very broad problem, so 
we will (for the time being) restrict to 
systems that have a ``good'' Lax representation.

\begin{defn}
An integrable system $(M,\chi)$ admits a ``Lax representation'' if 
it is isomorphic to a symplectic leaf of a 
meromorphic Hitchin system,
with $\chi$ the restriction of the Hitchin map.
\end{defn}

If the base curve has genus zero, this is essentially the same as the usual definition of a Lax representation.
Recall that the Hitchin map on 
moduli spaces of meromorphic Higgs bundles was shown to be proper in \cite{Nit-higgs}, and these moduli spaces  were shown to be integrable systems in 
a Poisson sense in \cite{Bot,Mar}.
It is this modular interpretation of integrable systems, as moduli spaces of Higgs bundles, initiated by Hitchin, that leads to the natural compactifications of the fibres.
More generally we will  consider meromorphic $G$-Higgs bundles and 
the extension allowing meromorphic Higgs fields on complex parahoric torsors.
However we don't want to consider arbitrary symplectic leaves of such Poisson moduli spaces; indeed the notion of Poisson integrable system does not imply every symplectic leaf is 
an integrable system in the sense of Defn. \ref{def: acihs}. 
Rather we will consider a special class of symplectic leaves that we will call ``good'' (see later below).
A key point is that many abstract integrable systems have more than one representation. 
This is essentially the story of different Lax pairs in the integrable systems literature. We view it as analogous to the fact that abstract Lie groups often have many faithful linear representations.

Two classes of examples of integrable systems are as follows 
(see e.g. the survey \cite{DonMar-short} for more background).

\subsection{Rational matrices}

Many classical examples of integrable systems fit into the framework of isospectral deformations of rational matrices.
Suppose $A(z)$ is an $n\times n$ matrix of rational functions on the Riemann sphere.
Then we can consider the coefficients of its characteristic polynomial
$$\chi = \det(A(z)-\la)$$
and this leads to the definition of the spectral curve.
One obtains symplectic varieties by considering rational matrices with the same orbits of principal parts
\beq\label{eq: mstar}
\cM^* = \{ \Phi= A(z)dz \st \text{orbits of polar parts fixed} \}/\GL_n(\IC)
\eeq
and for many such varieties the functions $\chi$ give an integrable system.
Jacobi's work on the geodesic flow on an ellipsoid 
fits into this set-up (see e.g. the exposition in \cite{DonMar-short}), 
and another early example  is due to Garnier \cite{garn1919}.
In these examples the spectral curve is a ramified covering of the Riemann sphere.
See also \cite{AHP, adler-vm-euc, Bea} for example.  

More explicitly consider the group $G_k:=
\GL_n(\IC[z]/z^k)$ and its Lie algebra $\g_k$.
Then using the trace pairing one can naturally view 
coadjoint orbits $\cO\subset \g_k^*$ as principal parts at $0$ of such matrices $A(z)dz$.
Repeating at each pole yields the identification
$$\cM^*\cong \cO_1\times \cdots \times \cO_m\spq 
\GL_n(\IC)$$
of $\cM^*$ with the symplectic quotient of the product of such orbits by the constant group $\GL_n(\IC)$, as in \cite{smid} \S2. In particular this gives $\cM^*$ a complex symplectic structure.

\subsection{Hitchin systems}

Given a Riemann surface $\Si$ of genus $g\ge 2$,
let $G=\GL_n(\IC)$, 
and consider the 
moduli space of stable rank $n$ vector bundles $\Bun_G$ (or more generally principal $G$-bundles for other complex reductive groups), and its cotangent bundle $T^*\Bun_G$.
Explicitly this is the space of pairs $(V,\Phi)$ where $V$ is a 
stable vector bundle and $\Phi\in \HH^0(\Si,\End(V)\otimes \Om^1)$ represents a cotangent vector to  $\Bun_G$ at $V$ (noting that the tangent space at $V$ is 
$\HH^1(\Si,\End(V))$).
More generally one can define a ``Higgs field'' on any vector bundle $V$ to be an element $\Phi\in \HH^0(\Si,\End(V)\otimes \Om^1)$, and then there is a stability condition for the pair $(V,\Phi)$, which is weaker than the condition for $V$ to be  stable.
This leads to a partial compactification of $T^*\Bun_G$, the moduli space $\MDol$ of stable Higgs bundles, whose points are isomorphism classes of stable pairs $(V,\Phi)$
(or S-equivalence classes of semistable pairs)\footnote{To simplify the presentation we will neglect to distinguish stability and polystability  throughout, referring the reader to the original articles for more precision.}.
This Higgs bundle moduli space is an integrable system, essentially taking the characteristic polynomial of $\Phi$, and its generalisation for other $G$ (see \cite{Hit-sbis}).

\begin{equation}	\label{diag: Higgs bundles}
\begin{array}{ccl}
 T^*\Bun_G & = &\{ (V,\Phi)\st \text{$V$ stable} \}/\text{isomorphism}  \\
\bigcap &&\\
\MDol & = &\{ (V,\Phi)\st \text{(poly)stable pair} \}/\text{isomorphism}  \\
\mapdown{\chi} &&\\
\IH &&
\end{array}
\end{equation}

\subsection{Nonabelian Hodge theory}

In a different paper the same year \cite{Hit-sde}, Hitchin also showed that the Higgs bundle moduli space had another beautiful property, that of being  hyperk\"ahler.
See \cite{hitchin-bbk} for a gentle introduction and many examples.
In brief this means the underlying manifold has a very special type of Riemannian metric 
with holonomy in the compact symplectic group, lying in the intersection of all the big classes of special holonomy groups on 
Berger's list:

$$
	
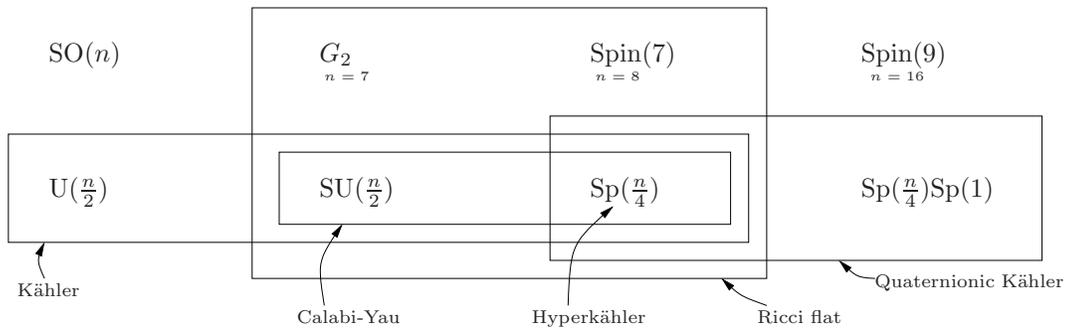
\begin{picture}(0,0)%
\includegraphics{holonomygroups.pstex}%
\end{picture}%
\setlength{\unitlength}{2486sp}%
\begingroup\makeatletter\ifx\SetFigFontNFSS\undefined%
\gdef\SetFigFontNFSS#1#2#3#4#5{%
  \reset@font\fontsize{#1}{#2pt}%
  \fontfamily{#3}\fontseries{#4}\fontshape{#5}%
  \selectfont}%
\fi\endgroup%
\begin{picture}(10329,3240)(1609,-6799)
\put(4726,-4111){\makebox(0,0)[lb]{\smash{{\SetFigFontNFSS{10}{12.0}{\rmdefault}{\mddefault}{\updefault}{$G_2$}%
}}}}
\put(4726,-5461){\makebox(0,0)[lb]{\smash{{\SetFigFontNFSS{10}{12.0}{\rmdefault}{\mddefault}{\updefault}{$\SU(\frac{n}{2})$}%
}}}}
\put(7426,-4111){\makebox(0,0)[lb]{\smash{{\SetFigFontNFSS{10}{12.0}{\rmdefault}{\mddefault}{\updefault}{$\Spin(7)$}%
}}}}
\put(7426,-5461){\makebox(0,0)[lb]{\smash{{\SetFigFontNFSS{10}{12.0}{\rmdefault}{\mddefault}{\updefault}{$\Sp(\frac{n}{4})$}%
}}}}
\put(2026,-5461){\makebox(0,0)[lb]{\smash{{\SetFigFontNFSS{10}{12.0}{\rmdefault}{\mddefault}{\updefault}{$\U(\frac{n}{2})$}%
}}}}
\put(2026,-4111){\makebox(0,0)[lb]{\smash{{\SetFigFontNFSS{10}{12.0}{\rmdefault}{\mddefault}{\updefault}{$\SO(n)$}%
}}}}
\put(4771,-4291){\makebox(0,0)[lb]{\smash{{\SetFigFontNFSS{5}{6.0}{\rmdefault}{\mddefault}{\updefault}{$n=7$}%
}}}}
\put(7471,-4291){\makebox(0,0)[lb]{\smash{{\SetFigFontNFSS{5}{6.0}{\rmdefault}{\mddefault}{\updefault}{$n=8$}%
}}}}
\put(6841,-6721){\makebox(0,0)[lb]{\smash{{\SetFigFontNFSS{7}{8.4}{\rmdefault}{\mddefault}{\updefault}{$\text{Hyperk\"ahler}$}%
}}}}
\put(9091,-6721){\makebox(0,0)[lb]{\smash{{\SetFigFontNFSS{7}{8.4}{\rmdefault}{\mddefault}{\updefault}{$\text{Ricci flat}$}%
}}}}
\put(4501,-6721){\makebox(0,0)[lb]{\smash{{\SetFigFontNFSS{7}{8.4}{\rmdefault}{\mddefault}{\updefault}{$\text{Calabi-Yau}$}%
}}}}
\put(10261,-6361){\makebox(0,0)[lb]{\smash{{\SetFigFontNFSS{7}{8.4}{\rmdefault}{\mddefault}{\updefault}{$\text{Quaternionic K\"ahler}$}%
}}}}
\put(1711,-6451){\makebox(0,0)[lb]{\smash{{\SetFigFontNFSS{7}{8.4}{\rmdefault}{\mddefault}{\updefault}{$\text{K\"ahler}$}%
}}}}
\put(10126,-4111){\makebox(0,0)[lb]{\smash{{\SetFigFontNFSS{10}{12.0}{\rmdefault}{\mddefault}{\updefault}{$\Spin(9)$}%
}}}}
\put(10126,-5461){\makebox(0,0)[lb]{\smash{{\SetFigFontNFSS{10}{12.0}{\rmdefault}{\mddefault}{\updefault}{$\Sp(\frac{n}{4})\Sp(1)$}%
}}}}
\put(10216,-4291){\makebox(0,0)[lb]{\smash{{\SetFigFontNFSS{5}{6.0}{\rmdefault}{\mddefault}{\updefault}{$n=16$}%
}}}}
\end{picture}%

$$

Hitchin established this by considering
the differential-geometric moduli space of
solutions of a certain system of nonlinear PDEs, 
the ``self-duality equations on a Riemann surface'',
now called the Hitchin equations.
This moduli space naturally has a \hk structure since it appears 
as the \hk reduction of an infinite dimensional \hk vector space;
indeed the Hitchin equations are the moment map for the action of the gauge group,
and so the set of gauge orbits of solutions is the \hk reduction. 
As explained in \cite{Hit-sde} \S1, 
these self-duality equations themselves are a dimensional reduction of the anti-self dual Yang-Mills/instanton  equations in dimension four, 
which can be viewed as the origin of their \hk nature.

Then in \cite{Hit-sde}, and more generally by Simpson \cite{Sim-hbls},
a bijective correspondence between stable Higgs bundles and solutions of Hitchin equations was established.
This correspondence is an instance of a general principle, 
interpreting stability conditions for algebro-geometric objects as the condition for existence 
of solutions of gauge-theoretic PDEs, the ``Hitchin--Kobayashi principle''.
There are lots of instances of this principle; 
 a simpler instance is the 
``Hitchin--Kobayashi correspondence for vector bundles on curves'', i.e.
the Narasimhan--Seshadri theorem \cite{NarSes}
(especially Donaldson's approach \cite{Don83}), 
relating 
stable algebraic vector bundles to unitary connections with constant central curvature (this curvature condition amounts to a nonlinear PDE, generalising the condition of being 
flat in the degree zero case).
Nowadays this result is subsumed as a special case of 
the Hitchin--Kobayashi correspondence for Higgs bundles, i.e. the case when the Higgs field is zero.

$$\boxed{\phantom{\bmx 1 \\ 1 \emx}
\MDol \quad \cong \quad \gM  = 
\{ \text{solutions of Hitchin's equations} \}/\text{isom.}
\phantom{\bmx 1 \\ 1 \emx}}$$
$$\text{\small Hitchin--Kobayashi correspondence for Higgs bundles, due to Hitchin and Simpson.}$$

This correspondence leads to a new perspective: 
\hk manifolds come equipped with a two-sphere of complex structures,
and in the case of $\gM$, this sphere is partitioned into three subsets: one point consisting of the Higgs bundle complex structure $\MDol$, its complex conjugate, and the remaining $\IC^*$ of complex structures, all isomorphic to the moduli space $\MDR$ of stable algebraic connections on vector bundles on $\Si$.

The proof of this statement amounts to establishing a different
instance of the 
Hitchin--Kobayashi principle, this time for 
algebraic connections on vector bundles.
Namely the stability condition on an algebraic connection turns out to be the condition for the existence of a twisted harmonic metric or ``harmonic bundle''. In turn the nonlinear equations for a harmonic bundle
are just another way to write out Hitchin's equations.
This is due to Donaldson \cite{Don87} (written as a companion to \cite{Hit-sde}) and more generally Corlette \cite{Cor88} (see also \cite{diederich-ohsawa85} for the case of $\SL_2(\IR)$).

$$\boxed{\phantom{\bmx 1 \\ 1 \emx}
\MDR \quad \cong \quad \gM  = 
\{ \text{solutions of Hitchin's equations} \}/\text{isom.}
\phantom{\bmx 1 \\ 1 \emx}}$$
$$\text{\small Hitchin--Kobayashi correspondence for connections, due to Corlette and Donaldson.}$$

Of course on a smooth compact complex algebraic curve the algebraic connections are the same as the holomorphic connections, and in turn to the complex flat $C^\infty$ connections, which are classified by their monodromy representations.
This gives a third natural algebraic structure on $\gM$ as the character variety, or Betti moduli space:

$$\boxed{\phantom{\bmx 1 \\ 1 \emx}
\MDR \quad \cong \quad \MB  = \Hom^\irr(\pi_1(\Si),G)/G 
\phantom{\bmx 1 \\ 1 \emx}
}$$
$$\text{\small Riemann--Hilbert correspondence for irreducible algebraic connections on $\Si$}$$

Note this does not change the complex structure, only the algebraic structure: the Riemann--Hilbert correspondence is a complex analytic isomorphism between two non-isomorphic non-compact algebraic varieties.
This third viewpoint gives the simplest description of the differentiable manifold underlying $\gM$:
\beq\label{eq: cv}
\Hom(\pi_1(\Si),G)/G
\cong 
\{ A_1,B_1,\ldots,A_g,B_g\in G\st
[A_1,B_1]\cdots[A_g,B_g]=1\}/G\eeq
where $[a,b]=aba^{-1}b^{-1}$, $g$ is the genus of $\Si$ and $G$ acts by conjugation, and 
$\rho\in \Hom(\pi_1(\Si),G)$ 
is irreducible if its image fixes no 
nontrivial proper subspace of $\IC^n$.

Thus in summary there is a rich picture of a space with three natural 
algebraic structures, where the three maps on the right are isomorphisms.

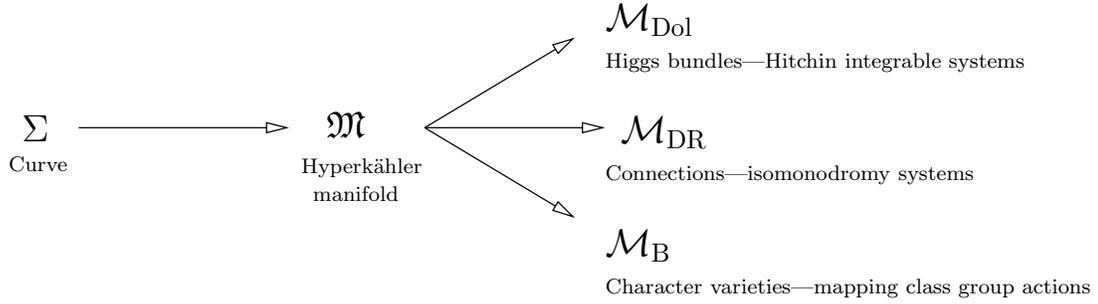
\begin{figure}[h]
\input{scheme.pstex_t}
\end{figure}

The terminology {\em Dolbeault, De\! Rham, Betti} for these algebraic
varieties comes from the viewpoint of nonabelian cohomology, and these spaces are three realisations of (the coarse moduli space underlying) the first cohomology $\HH^1(\Si,G)$ with $G$ coefficients. See Simpson \cite{Sim-hbls, Sim94ab, simp-hfnc} for more on this motivic viewpoint and for a general construction of the algebraic moduli spaces.
The resulting comparison isomorphism
$$\MDol\cong \MDR$$ is 
(thus) called the nonabelian Hodge correspondence.
Note that it is the composition of two {\em different} Hitchin--Kobayashi correspondences (proved by different people).
The other comparison isomorphism $\MDR\cong \MB$ is Riemann--Hilbert, and is used to define the isomorphism $\gM\cong \MB$ in the big diagram above.

$$\boxed{ \phantom{\bmx 1 \\ 1 \\ 1 \emx} \MDol
\overcong{\text{Nonabelian Hodge}}  
\MDR 
\overcong{\text{Riemann--Hilbert}} 
\MB
\phantom{\bmx 1 \\ 1 \emx}}$$

\subsection{Isomonodromy and nonabelian Gauss--Manin connections}\label{sn: imd}

The Higgs bundle moduli spaces $\MDol$ have natural flows on them---the Hamiltonian flows of the integrable systems: Each component of the Hitchin map is a function on the holomorphic symplectic manifold $\MDol$, and so has an associated Hamiltonian vector field (when written in local coordinates these amount to differential equations)---the corresponding flow is a ``straight line flow'' in each of the abelian fibres.

On the other hand the moduli spaces of connections $\MDR$ are the arena for a different family of flows.
More precisely if we have a smooth family of compact Riemann surfaces 
$\underline\Si\to \IB$ over a base $\IB$, with fibres $\Si_b$ for $b\in \IB$, 
then we can consider the corresponding family of  moduli spaces of connections
$\bcM\to \IB$ whose fibre over $b\in \IB$  is $\MDR(\Si_b)$.
The flows are on the total space $\bcM$:

\begin{thm}[\cite{Sim94ab}]
The fibre bundle $\bcM\to \IB$ has a natural complete flat (algebraic) Ehresmann connection on it, the nonabelian Gauss--Manin connection.
\end{thm}

If written explicitly in local coordinates this amounts to a
system of nonlinear differential equations. 
If $\dim(\IB)=1$ it amounts 
to a line field on $\bcM$ transverse to the fibres.
Integrating this connection gives a natural analytic way to identify fibres $\MDR(\Si_a)\cong \MDR(\Si_b)$, for any path in $\IB$ from $a$ to $b$.
The same identification can be obtained via the Riemann--Hilbert correspondence, passing to the Betti spaces, and then identifying them by dragging loops around, and keeping the monodromy representation constant---leading to the familiar mapping class group actions on the character varieties.
This fits together Simpson's viewpoint on nonabelian Gauss--Manin connections and the more classical viewpoint of isomonodromic (monodromy preserving) deformations of linear connections (as explained in a more leisurely way in \cite{smid}).
We will see more examples of this below, and discuss how it generalises
to the case of meromorphic connections, leading to the notion of wild Riemann surfaces and wild mapping class groups.

The upshot is that the nonabelian Hodge package encodes two types of nonlinear differential equations: 
the integrable systems and the isomonodromy systems.
Both of these types of differential equations have a Lax problem:
1) to find a rational matrix (or more generally a Higgs bundle) whose isospectral deformations are controlled by the given nonlinear equation, or 2) to find a linear differential system/connection whose isomonodromic deformations are controlled by the given nonlinear equation.
And in many cases there is more than one distinct Lax representation of the nonlinear equation (typically an infinite number---see e.g. \cite{slims} \S11.3).

Thus it makes sense to try to classify the ``representations'' of the whole nonabelian Hodge package, and not just the Higgs bundle  or connection moduli space. 

\subsection{Extending the nonabelian Hodge package}

Of course to really fit all this together we need to see that the
classical integrable systems on the spaces $\cM^*$ of rational matrices
can indeed be viewed on the same footing as the sophisticated nonabelian Hodge set-up.

An $n\times n$  
rational matrix $A(z)$ (representing a point of $\cM^*$) 
yields a matrix $A(z)dz$ of meromorphic one-forms.
We can perfectly well view this as a meromorphic Higgs field on the trivial rank $n$ vector bundle on $\IP^1$.
We will write $\cM^*=\MDol^*$ when we think of it in this way.
Similarly $A$ determines a meromorphic connection $d-A(z)dz$ on the trivial rank $n$ vector bundle on $\IP^1$.
We will write $\cM^*=\MDR^*$ when we think of it in this way.
Thus we have lots of holomorphic symplectic moduli spaces of connections and Higgs bundles with arbitrary order poles.

Thus we can naively ask if there are extensions: 
1) of the Riemann--Hilbert correspondence to include $\MDR^*$, or 2) of nonabelian Hodge yielding $\MDol^*\cong \MDR^*$.

The answer is {\em no}, for several reasons.
We will sketch what can be done however, and how to adjust the question (see  \cite{ihptalk} for a more detailed review).
The first issue is that it is too stringent to insist the underlying vector bundle is holomorphically trivial, even in the case of $\IP^1$, when such bundles are generic amongst topologically trivial bundles.
Rather one should just fix the topological type.

Let $\Si$ be a compact Riemann surface and $\ba=\{a_i\}\subset \Si$ a finite set of marked points. Fix integers $n$ (the rank) and $k_i\in\IZ_{\ge 1}$ for each marked point, and let $D=\sum k_i(a_i)$ be the resulting positive divisor.
Nitsure \cite{Nit-higgs} constructs algebraically a moduli space $\MDol$ of meromorphic 
Higgs bundles with poles on $D$, and shows the corresponding Hitchin map is proper. 
Similarly moduli spaces of meromorphic connections may be constructed 
(\cite{Nit-log}, if each $k_i=1$, and one may use \cite{Sim94ab} in general).
A Poisson structure on $\MDol$ was constructed in \cite{Bot, Mar}, showing it is an integrable system in a Poisson sense.
The symplectic leaves $\MDol(\bcO)\subset \MDol$ are obtained by fixing the $G_k$ orbits $\cO_i$ at each pole, as in the case of $\cM^*$. 
In the case $\Si=\IP^1$ we get a map $\MDol^*\to \MDol(\bcO)$ onto a Zariski open subset.
A more reasonable question to ask is if one can extend the nonabelian Hodge correspondence to $\MDol(\bcO)$: are these spaces \hk (becoming moduli spaces of meromorphic connections $\MDR(\bcO)$ in another complex structure in the \hk family)? Is there a Riemann--Hilbert correspondence for such connections?

In this degree of generality, for arbitrary orbits $\bcO$, the answer is not known. However if we restrict the orbits a little bit, then we can proceed, as follows.

Fix a positive integer $k\ge 1$ and a coadjoint 
orbit $\cO\subset \g_k^*$, in the dual of the 
Lie algebra of $G_k=G(\IC[z]/z^k)$, where $G=\GL_n$. 
Let $\lt\subset \g$ be a fixed Cartan subalgebra (such as the diagonal matrices).
As in \cite{smid}
we view elements of $\g_k^*$ 
as principal parts of meromorphic connections/Higgs fields at a pole of order at most $k$ (with local coordinate $z$).

\begin{defn}
An orbit $\cO\subset \g_k^*$ is ``very good'' if it contains an element of the form $$dQ + \La\frac{dz}{z}$$
for some $\La\in \g$ and element $Q = \sum_1^{k-1}A_i/z^i\in \lt\flp z \frp/\lt\flb z \frb$.
\end{defn}

The diagonal element $Q$ is the ``irregular type''.
Lets fix an irregular type $Q_i$ at each marked point $a_i$, with 
pole of order $k_i-1$ 
(a coordinate independent approach is possible \cite{gbs}). 
Clearly not every orbit is very good, but most of them are, for example if the leading term is regular semisimple.
Given an irregular type $Q$ let $H\subset G$ be its centraliser (the subgroup commuting with each coefficient $A_i$).
Using the $G_k$ action we can (and will) assume $\La\in \lh$ is in the Lie algebra of $H$. 
Fix $\La_i\in \lh_i$ for each marked point $a_i$, and let 
$\cO_i\subset \g_{k_i}^*$ be the corresponding orbit, 
associated to $Q_i,\La_i$.

The simplest case to state is that with each formal residue $\La_i$ zero. 
Let $\cO_i'$ be the orbit associated to $-Q_i/2$ (with $\La_i=0$).
Thus we have orbits $\bcO= \{\cO_i\}$  and $\bcO'= \{\cO_i'\}$, and thus spaces of meromorphic connections $\MDR(\bcO)$ and Higgs bundles $\MDol(\bcO')$, with fixed principal parts.

\begin{thm}[\cite{wnabh, Sab99}]
There is a moduli space $\gM$ of solutions to Hitchin's equations on 
$\Si\setminus \ba$ which is a \hk manifold,
isomorphic to $\MDR(\bcO)$ in one complex structure
and to $\MDol(\bcO')$ in another.
\end{thm}

In brief \cite{Sab99} establishes the Hitchin--Kobayashi correspondence for meromorphic connections, and \cite{wnabh}
constructs the moduli spaces $\gM$ and establishes the Hitchin--Kobayashi correspondence for meromorphic Higgs bundles, thus establishing the wild nonabelian Hodge correspondence on curves.
At the other extreme is the case when each $Q_i=0$, so $k_i=1$ and every orbit is very good.
The nonabelian Hodge correspondence was extended to this tame/logarithmic 
case earlier by Simpson \cite{Sim-hboncc}.
A key subtlety here is that {\em one needs to incorporate parabolic structures at the poles to get a complete correspondence}.
In particular it becomes clear that the weights appear on the same footing as the real/imaginary parts of the eigenvalues of the residues
(and in turn the resulting triple is best thought of in terms of imaginary quaternions).
More explicitly the parabolic weight on the Higgs side gives the real part of the eigenvalues of $\La$ on the connection side (for example).
Simpson phrases this in terms of filtered bundles, and that leads to the notion of parahoric bundle, that works equally well with other structure groups 
(see \cite{logahoric}).
In the notation of \cite{logahoric} the full table that Simpson found in 
\cite{Sim-hboncc} is as follows:

  \begin{center}
  \begin{tabular}{| c || c | c | c |}
    \hline
   \rule{0pt}{2.4ex}      & \text{Dolbeault} & \text{De\! Rham} & \text{Betti}\\ \hline
    \rule{0pt}{2.4ex}\text{weights}\,$\in \lt_\IR$ & 
$-\tau$ & $\th$ & $\phi=\tau+\th$ \\ \hline
  \rule{0pt}{2.5ex}\text{ eigenvalues}\,$\in \lt_\IC, \lt_\IC, \text{T}(\IC)$ & 
$-(\phi+\si)/2$ & $\tau+\si$ & $\exp(2\pi i (\tau+\si))$ \\
    \hline
  \end{tabular}
  \end{center}

The same rotation of the weights/eigenvalues (of the formal residues $\La_i$) persists in the wild case \cite{wnabh}, superposed onto the change in irregular type 
$-Q_i/2\leftrightsquigarrow Q_i$ already described.

For example in the De\! Rham column this means: a weight 
$\th\in \lt_\IR$ determines a parahoric subalgebra 
$\wp_\th\subset \g\flp z \frp$ and 
we define a connection to be ``$\th$-logahoric'' if
it is locally of the form  
$\Th(z)\frac{dz}{z}$ with $\Th\in \wp_\th$.
In brief $\th$ determines a grading on
$\g\flp z \frp$, and $\wp_\th$ is
the non-negative piece.
For example if $\th=0$ then $\wp_\th = \g\flb z \frb$ 
and so we recover the notion of logarithmic connection.
More generally if the components of $\th$ are in
the interval $[0,1)$ this is a 
logarithmic connection on a parabolic vector bundle (with the residue preserving the flag in the fibre)---the 
elements $\tau,\si$ are then the real and imaginary parts of the eigenvalues of $\La=\Th(0)$.
Then we look at the ``very good''  connections, which are locally of the form 
\beq
dQ+\Th(z)\frac{dz}{z}, \qquad 
\Th\in \wp_\th\subset \g\flp z \frp,
\eeq
i.e. of the form  ``$dQ+\th\text{-logahoric}$'', with $\th\in \lt_\IR, Q\in\lt\flp z \frp$.
Similarly for $\MDol$.

\begin{rmk}
For $G=\GL_n$ one can always act with the loop group 
to reduce to the parabolic case, but this is not true for general reductive groups.
Nonetheless the above definition of 
$\th\text{-logahoric}$ (from \cite{logahoric}) makes sense, and we can define the very good connections in the same way.
More generally one can define a meromorphic connection on a parahoric bundle to be ``good'', if locally at each pole there is a cyclic cover $z=t^r$ 
such that the connection becomes very good after pullback. 
In the case of $\GL_n$ such twisted/ramified connections were already considered in \cite{Sab99} 
(see also e.g \cite{bertola-mo, bremer-sage-IMD, inaba-twisted}), and the analysis in \cite{wnabh} still works.
For other groups we conjecture this is the right class of connections to look at, from the viewpoint of nonabelian Hodge theory and Riemann--Hilbert\footnote{One 
can also ask: are the ``good'' meromorphic Higgs bundle moduli spaces exactly the ones which are integrable systems in the sense of Defn. \ref{def: acihs}? We leave this as a question rather than a conjecture as we've not looked into it. This is compatible with 
\cite{BKV} though, which came to light during conference (although they use a more stringent definition).}.
\footnote{The terminology for parahoric bundles is similar to that for parabolic bundles: 
A {\em quasi-parahoric bundle} is a torsor for a parahoric group scheme $\cG\to \Si$, as in 
\cite{pap-rap2}. Locally such a group scheme amounts to a parahoric subgroup $\cP$ of a formal loop group $G\flp z \frp$ (or a twisted loop group), as in 
\cite{BrTits-I, pap-rap}, and so by definition such a torsor is locally isomorphic to such a subgroup.
In turn a {\em weighted parahoric subgroup} of the loop group is a {\em point} of the corresponding Bruhat-Tits building, in the facette corresponding to $\cP$ (see \cite{logahoric} Defn. 1 p.46).
Up to conjugacy by the loop group this intrinsic definition of weight reduces to the naive notion of weight $\th\in \lt_\IR$---indeed the building is built out of apartments $\lt_\IR$.
Finally a {\em parahoric bundle} is a quasi-parahoric bundle plus a choice of a weight at each point.
More recent references include \cite{BS-pub, B-GP-MR, heinloth-stability}.
}

\end{rmk}

\subsection{Nonabelian Hodge spaces}

Thus, rather than just classifying integrable systems, 
we could try to 
classify the richer geometric objects occurring in this nonabelian Hodge story. 
For this it is convenient to make the following definition. 

\begin{defn}\label{def: nabhs}
A ``nonabelian Hodge space'' is a \hk manifold $\gM$ with three preferred algebraic structures $\MDol,\MDR,\MB$, such that $\MDol$ is an algebraic integrable system in the sense of Definition \ref{def: acihs}, and is a symplectic leaf of a meromorphic Hitchin system.
\end{defn}
In the first instance we will focus on 
the {\em complete} nonabelian Hodge spaces, 
i.e. those whose \hk metric is complete (see \cite{wnabh} for sufficient conditions to ensure this).
More generally we will allow $\gM$ to have certain 
(orbifold) singularities 
(as occur even for stable $G$-Higgs bundles, once we move away from $\GL_n(\IC)$).
From the preceding discussion (and existing results on the irregular Riemann--Hilbert correspondence 
\cite{sibuya77, Mal79, deligne78, jurkat78, BV89, malg-book, MR91, L-R94, DMR-ci}) 
we know there are lots of examples of interest.
As usual by \hk rotation the integrable system on $\MDol$ yields 
a special Lagrangian fibration
on $\MDR\cong \MB$.

\begin{rmk}
Note that here we are tacitly restricting to connections/gauge theory/Higgs bundles on Riemann surfaces/smooth complex algebraic curves.
Our viewpoint on wild nonabelian Hodge theory is to use it to produce new
moduli spaces as output, and then study them and the nonlinear differential equations that live on them. 
Note that lots of work has been done recently,
with quite different motivation, 
to extend the nonabelian Hodge correspondence  to higher dimensional varieties
(cf. \cite{Cor88, Sim-hbls} in the compact case, 
\cite{Biq97, jostzuo97, mochizuki-tame} in the noncompact case with tame singularities, 
and \cite{mochizuki-wild} for the Hitchin--Kobayashi
correspondence for irregular connections on quasi-projective varieties).
Note that it is not at all clear if any new nonabelian Hodge spaces occur in the higher-dimensional set-up, so we leave this as 
a provocative open problem
(motivated by \cite{Sim04} p.2):

\noindent {\bf Problem.}\ 
Find an example of a nonabelian Hodge space arising as a component of a moduli space of connections on a smooth quasi-projective variety, that is not isomorphic to one arising from a curve.
\end{rmk}

\section{Non-perturbative symplectic manifolds}

To get a feel for these spaces we will describe some of the underlying holomorphic symplectic manifolds, from the Betti perspective, 
which is often the most concrete description.
 The \hk approach of \cite{wnabh} strengthens the  earlier complex symplectic construction \cite{thesis, smid}.
Fix $\Si, D$ and irregular types $\bQ=\{Q_i\}$ as above.
Write $\bSi = (\Si,D,\bQ)$ for this triple, an irregular curve/wild Riemann surface.
The simplest examples are when $\Si=\IP^1$, and then we can view 
$\cM^*$ as an approximation to the full moduli space.
First suppose we are in the tame case $Q_i=0, k_i=1$.

The tame character varieties generalising \eqref{eq: cv} 
are as follows:
Choose a conjugacy class $\cC_i\subset G$ for each marked point
and write $\cC=\{\cC_i\}$.
Then $\MB(\bSi,\cC)=\Hom_\cC(\pi_1(\Si^\circ),G)/G$
where $\Hom_\cC(\pi_1(\Si^\circ),G)$ is isomorphic to
\beq \label{eq: tamecv}
\{ \bA,\bB,\bM\in G^{g}\times G^g \times \cC \ \st\ 
[A_1,B_1]\cdots[A_g,B_g]M_1\cdots M_m=1\}.
\eeq
They fit into Deligne's Riemann--Hilbert correspondence \cite{Del70} (on $\Si^\circ)$, and also into Simpson's tame nonabelian Hodge correspondence 
\cite{Sim-hboncc} (upon taking the Betti weights $\phi$ zero). As described in op. cit., in general one should consider filtered local systems---the resulting character varieties will not always be affine or quasi-affine (this amounts to replacing the classes $\cC_i$ by weighted conjugacy classes $\wh \cC_i$ from \cite{logahoric}).

The simplest nontrivial example is rank two with four marked points:

{\bf 1) The Fricke--Klein--Vogt surfaces.}
Let $G=\SL_2(\IC)$, and $\bSi=(\IP^1,\ba)$ be the sphere with marked points $\ba=(a_1,\ldots,a_4)=(0,t,1,\infty)$.
Let  $\Si^\circ = \IP^1\setminus\ba$ and
choose regular semi-simple conjugacy classes $\cC_i\subset G, i=1,\ldots,4$.
The full character variety $\Hom(\pi_1(\Si^\circ),G)/G$ is a six-dimensional Poisson variety and its generic symplectic leaves are of the form
\beq\label{eq: fkvredn}
(\cC_1\fus{} \cdots \fus{} \cC_4)\spq G = 
\{ (M_1,\ldots,M_4)\in G^4\st M_i\in \cC_i, M_1\cdots M_4=1\}/G.
\eeq
These are affine complex surfaces, given 
(\cite{vogt} eq. (11), \cite{FrickeKleinI} p.366, \cite{Magnus}) 
by an  equation of the form:
\beq\label{eq: fkv}
xyz+x^2+y^2+z^2 + ax+by+cz=d
\eeq
for constants $a,b,c,d\in \IC$ determined by the 
eigenvalues of the $\cC_i$.
The quotient \eqref{eq: fkvredn} is a quasi-Hamiltonian or {\em multiplicative symplectic quotient}, involving group valued moment maps as in \cite{AMM}.
The corresponding additive symplectic quotient
is one of the spaces $\MDR^*$: choose four 
coadjoint orbits $\cO_i\subset \g\cong \g^*$.
Then
\beq\label{eq: addfkv}
(\cO_1\times \cdots \times \cO_4)\spq G = 
\{ (A_1,\ldots,A_4)\in \g^4\st A_i\in \cO_i, \sum A_i=0\}/G
\eeq
which is of the form \eqref{eq: mstar} with 
$A = \sum_1^3\frac{A_i}{z-a_i}$. %
Then we can take the monodromy of the connection $d-Adz$ 
to get a holomorphic map 
\beq\label{eq: nua}
\nu_\ba: (\cO_1\times \cdots \times \cO_4)\spq G \ \to\   
(\cC_1\fus{} \cdots \fus{} \cC_4)\spq G
\eeq
from the additive to the multiplicative symplectic quotient 
(if $\cC_i=\exp(2 \pi \sqrt{-1}\cO_i)$, and if none of the 
residual eigenvalues differ by a nonzero integer).

\begin{thm}[Hitchin \cite{Hit95long}]\label{thm: nusptame}
For any choice of points $\ba\subset \IP^1$ (and generators of the fundamental group) the transcendental 
map $\nu_\ba$ is a holomorphic symplectic map.
\end{thm}

Thus we see that the Atiyah--Bott/Goldman symplectic structure on the character variety has the somewhat magical property of matching that on the additive space, for {\em any} choice of pole configuration.
This property holds in much more generality, 
even for irregular singular connections/wild character varieties, and even when the deformation space (of curve with marked points) is generalised to an irregular curve.

If we vary the curve-with-marked-points $\bSi$, we are really just moving $t\in \IB=\IP^1\setminus\{0,1,\infty\}$. 
Then, as in \S\ref{sn: imd},  there is an isomonodromy connection  on the bundle 
$\bcM\to \IB$ of De\! Rham spaces, 
with fibres of complex dimension two.
Written explicitly this isomonodromy connection
becomes a second order differential equation---in this case the Painlev\'e VI equation (see e.g. \cite{srops}).
Thus we have a link between \hk four-manifolds and Painlev\'e equations.
The next (irregular) example corresponds to taking the Painleve II equation:

{\bf 2) The Flaschka--Newell surfaces}---multiplicative Eguchi--Hanson spaces.

Recall the Eguchi--Hanson space \cite{EH78} was the 
first nontrivial example of a complete \hk manifold, and is $T^*\IP^1$ in one complex structure and an affine $\SL_2(\IC)$ coadjoint orbit in its generic complex structure.
Somewhat improbably, it occurs as one of the spaces $\cM^*$ as follows.
Suppose $G=\SL_2(\IC)$ and $\bSi=(\IP^1,\infty,Q)$
with just one marked point, with irregular type 
$Q=A_3z^3+A_2z^2+A_1z$, 
having a pole of order $3$ at $\infty$, with $A_3$ regular.
The resulting connections have a pole of order $4$ at $\infty$.
Fix nonzero $\La\in \lt$.
The corresponding space $\cM^*$ has complex dimension two.
Indeed, if $w=1/z$, 
the group $G_4:=G(\IC[w]/w^4)$ of jets at $z=\infty$
has Lie algebra 
$\g_4=\g(\IC[w]/w^4)=\{X= \sum_0^3 X_iw^i\st X_i\in \g\}$ and the dual of this can be identified with 
$$\g_4^* = \{B=\sum_1^4 B_i\frac{dw}{w^i}\st B_i\in \g\}$$
via the pairing $\langle X,B \rangle = \res \tr(XB)$.
As usual we identify $dQ+\La dw/w$ as a point of $\g_4^*$, and 
let $\cO\subset \g_4^*$ be its coadjoint orbit,
a holomorphic symplectic manifold of dimension $8$.
This has a Hamiltonian action of $G$ via the inclusion 
$G\subset G_4$ and the coadjoint action, with moment map given by the residue. Note the elements with zero residue are precisely those that extend holomorphically to the finite $z$-plane, and so
the symplectic quotient is the space $\cM^*$:
\beq \label{eq: irraddspq}
\cM^* \cong \cO\spq G.
\eeq
This has dimension $8-2\dim(G)=2$.

\begin{lem}[\cite{quad} ex.3, \cite{rsode} Apx. C, \cite{hi-ya-nslcase}]
$\cM^*$ is isomorphic to the Eguchi--Hanson space in its generic complex structure as a complex symplectic manifold.
\end{lem}

The theory of Stokes data, properly interpreted, 
then gives us a complex surface $\MB$, the wild character variety, 
and a holomorphic map $\nu_Q:\cM^*\to \MB$,
generalising the maps taking the monodromy representation in the tame case. 
The space $\MB$  can be equipped with a holomorphic 
symplectic structure, such that one again has the 
analogous magical property.
In this example it is thus the multiplicative analogue of the 
Eguchi--Hanson space.
The underlying algebraic surface $\MB$
was written down by Flaschka--Newell \cite{FN80} (3.24) (in fact using a different Lax pair for Painlev\'e II) as the affine surface
\beq\label{eq: fn}
xyz + x+y+z = d
\eeq
for a constant $d\in \IC$.

To explain the extension of Thm. \ref{thm: nusptame}
we need to discuss more the notion of irregular curve generalising the underlying curve with marked points in the tame case above.
Indeed in this example the underlying curve with marked point 
$(\IP^1,\infty)$ has no moduli, but the  irregular curve
$(\IP^1,(\infty,Q))$ lives in a one-dimensional moduli space:
There are $3$ parameters in $Q$, but we can act with the $2$ dimensional group of Mobius transformations 
fixing $\infty$, leaving one parameter, an irregular analogue of the cross-ratio of the four points in the tame case.
On the other hand the constant $d$ in \eqref{eq: fn} is determined by the parameter $\La$ (the exponent of formal monodromy), similarly to the constants $a,b,c,d$ in the tame case.
The analogue of Thm \ref{thm: nusptame}, will thus arise when 
we keep $\La,d$ fixed, and allow $Q$ to vary:

\begin{thm}[\cite{smid} Thm 6.1] \label{thm: irhmap}
The transcendental map $\nu_Q : \cO\spq G \to \MB$ taking the connections in $\cM^*\cong \cO\spq G$ 
 to their Stokes data, is a symplectic map
for any choice irregular type $Q$ (and other discrete choices). 
\end{thm}

This holds quite generally and thus we have a larger class of symplectic manifolds with similar magical properties to the tame case (they are the symplectic manifolds in the title of \cite{smid}---the name ``wild character variety'' is more recent).

\section{New quasi-Hamiltonian spaces}

In the tame case \eqref{eq: nua} we had a map
$\nu_\ba:\cO_1\times\cdots\times \cO_4\spq G \to 
\cC_1\fus{} \cdots \fus{} \cC_4\spq G$
so the Betti spaces clearly {\em looked like} a multiplicative version of the additive side.
A similar picture holds in the irregular case, for example involving  the multiplicative analogue of the symplectic description 
$\cM^*\cong \cO\spq G$ in \eqref{eq: irraddspq}.
Thus we will here describe some simple spaces of Stokes data and refer to the literature for more general examples.
See \cite{p12} or \cite{even-euler} for a more comprehensive review.

First we rephrase the additive side, showing how it decouples.
Suppose $G=\GL_n$ and $\cO\subset \g_k^*$ is a very good orbit, containing  $dQ+\La dz/z$ say.
The orbit $\cO$ may be decoupled as follows (cf.
\cite{smid} \S2).
Let $B_{k}\subset G_{k}$ be the kernel of the evaluation map 
$G_{k}\onto G$.
Then $dQ$ may be viewed as a point of the dual of the Lie algebra of
$B_{k}$.
Let $\cO_B\subset \Lie(B_{k})^*$ be its coadjoint orbit.
Let $\wt \cO= \cO_B\times T^*G$.
We call this space $\wt \cO$ the ``extended orbit''; it is a
Hamiltonian $G\times H$-space.
It arises by allowing the formal residue $\La$ to vary and adding a compatible framing---the term ``extended'' is
by analogy with the set-up of
Jeffrey \cite{Jef94}.
The orbit $\cO$ arises as the reduction
$\wt \cO \spq_\La H$ at the value $\La$ of the moment map.
On the other hand $G$ only acts on $T^*G$ and 
the reduction $\wt \cO\spq G$
is isomorphic to $\cO_B$.
For example this implies:
\begin{cor}[\cite{smid}]
In the set-up of \eqref{eq: irraddspq}
the space $\cM^*\cong \cO\spq G$ 
is isomorphic to the reduction $\cO_B\spqa{\La} H$
of $\cO_B$ at the value $\La$ of the moment map for $H$.
\end{cor}

In particular, as noted in \cite{smid}, 
a theorem of Vergne \cite{Vergne72} implies 
 $\cO_B$ has global Darboux coordinates, 
 which leads to the fact that $\cM^*$ does not change under deformations of $Q$
 (more precisely under ``admissible'' deformations, defined in \cite{gbs}).

Now in the multiplicative case we will describe three spaces 
$\cA,\cB,\cC$ which are the multiplicative analogues of 
the spaces $\wt \cO, \cO_B, \cO$ respectively.
A key point is that: {\em $\cC$ is not a conjugacy class of $G_k$
if $k>1$}.
They have quasi-Hamiltonian actions of $G\times H, H, G$ respectively, and $\MB$ will have three descriptions:
$$\MB
\quad\cong\quad
 \cC\spq G 
           \quad\cong\quad   
 \cB\spqa{q}H 
                 \quad \cong\quad   
G\sqp\cA\spqa{q}H$$
where $q=\exp(2\pi i \La)\in H$, analogous to
$$\cM^*
\quad\cong\quad
 \cO\spq G 
           \quad\cong\quad   
 \cO_B\spqa{\La}H 
                 \quad \cong\quad   
G\sqp\wt \cO\spqa{\La}H.$$

The first step is:
\begin{thm}[\cite{saqh-ap}]
Suppose $Q=\sum_1^rA_i/z^i$ has regular 
semisimple leading term $A_r$, where $r=k-1$.
Let $U_\pm$ be the unipotent radicals of a pair of opposite Borels 
$B_\pm\subset G$.
Then the ``fission space'' 
$$\cA=\cA(Q) = G\times (U_+\times U_-)^r\times H$$
is an algebraic quasi-Hamiltonian $G\times H$ space.
\end{thm}

For example if $G=\GL_2$ or $\SL_2$ 
$$
U_+ = \bmx 1 & * \\ 0 & 1 \emx,  \
U_- = \bmx 1 & 0 \\ * & 1 \emx, \
H = \bmx * & 0 \\ 0 & * \emx \subset G.
$$
In fact \cite{saqh-ap} proves this for arbitrary complex reductive groups $G$. 
The space $\cA$ is denoted $\wt \cC/L$ in \cite{saqh-ap} Rmk 4 p.6.
The simplest fully nonabelian extension 
(with $B_\pm$ replaced by arbitrary opposite parabolics, and $H$ by their common Levi subgroup) appears in \cite{fission}. 
More generally see \cite{gbs} Thm 7.6 for the fission spaces $\cA(Q)$ for arbitrary irregular types $Q$, and 
\cite{twcv} for the twisted case.

The spaces $\cB=\cA\spq G$ and $\cC=\cA\spq_{\!q\,} H$ follow from this.
For example 
$$\cB  = \{ (\bS,h)\in(U_+\times U_-)^r \times H\st
 hS_{2r}\cdots S_2S_1 =1\}$$
which is a quasi-Hamiltonian $H$-space with moment map $h^{-1}$.
Here $\bS = (S_{1},\ldots,S_{2r})$ with $S_{odd/even}\in U_{+/-}$ respectively.

For example it is now a simple exercise, in the case 
$r=3,G=\SL_2(\IC)\supset H\cong \IC^*$, to compute the quotient 
$\MB=\cB\spq_{\!q}\, H$ and obtain the Flaschka--Newell surface.

The map $\nu_Q$ of Thm. \ref{thm: irhmap} now goes from 
$\cO\spq G$ to $\cC\spq G$, and similarly for any number of poles on $\IP^1$, with $\cO_i=\cO_i(Q_i,\La_i)$ and $\cC_i=\cC_i(Q_i,\La_i)$:
$$\nu_{\ba,\bQ}:\left(\cO_1\times\cdots\times \cO_m\right)\spq G 
\ \to \ 
\left(\cC_1\fus{} \cdots \fus{} \cC_m\right)\spq G.$$ 
The wild character varieties of any irregular curve 
$\bSi=(\Si,D,\bQ)$ %
are similar:
$$\MB(\bSi,\bLa) \cong 
\left(\ID^{\fus{}g}\fus{} \cC_1\fus{}\cdots \fus \cC_m\right)\spq G$$
where $g$ is the genus of $\Si$ and 
$\ID=G^2$ (which is a quasi-Hamiltonian $G$-space with moment map $aba^{-1}b^{-1}$).
For generic elements $\La_i$ 
the spaces $\MB(\bSi,\bLa)$ are smooth symplectic algebraic varieties 
and a formula for their dimension is in  \cite{gbs} Rmk. 9.12.
A more intrinsic approach, involving the space of Stokes representations $\Hom_\IS(\Pi,G)$, generalising the familiar space of fundamental group representations,
appears in \cite{gbs} and is reviewed in \cite{p12}.

The theory of irregular isomonodromy \cite{garn1912, FN80, JMU81}
lies behind this notion of irregular curve/wild Riemann surface 
(showing that the irregular type is on the same footing as the moduli of the curve),
although the recent extensions \cite{slims, gbs}
go beyond the classical theory by allowing arbitrary irregular types (e.g. the leading coefficient may have repeated eigenvalues).
The underlying idea is justified by the following.

\begin{thm}[\cite{smid, saqh-ap, gbs}]
Suppose $\underline{\bSi}\to \IB$ is an admissible family of irregular curves over a base $\IB$.
Then the corresponding bundle of wild character varieties 
$\bcM\to \IB$ is a local system of Poisson varieties.
In particular it carries a natural complete flat Poisson Ehresmann connection (the irregular isomonodromy connection). 
\end{thm}

In particular this implies the fundamental group of $\IB$ acts on the 
the fibres $\MB(\bSi,\bLa)$ by algebraic Poisson automorphisms.
For example many $G$-braid group actions occur this way \cite{bafi}.
Considering all such families  leads to the moduli stack of admissible deformations of a given wild Riemann surface, and in turn to the wild mapping class group action on the wild character varieties.

\section{More examples}

Thus the aim is to construct a large table
of nonabelian Hodge spaces and find which are isomorphic (or isogenous, 
deformations of each other, etc).
The first step can be initiated by 
glancing at papers/books on integrable systems or isomonodromy, and might begin as follows:

\renewcommand{\arraystretch}{1.2}
\begin{center}
  \begin{tabular}{| c | c | c | c |}
 \hline
     Rational map     & 
    \shortstack{\  \\ integrable \\system}
         & 
    \shortstack{isomonodromy\\system} &
    \shortstack{\ \\(wild) character \\ variety}
\\
 $\Phi$ & $\MDol$ & $\MDR$ & $\MB$
 \\ \hline\hline 
 $(A+Bz)\frac{dz}{z}$ & Mi\v s\v cenko-Fomenko & 
$\phantom{\begin{smallmatrix} \! \\ \! \\ \! \\ \!  \end{smallmatrix}}
\underset{\text{(Jimbo-Miwa-M\^{o}ri-Sato)}}{\text{Dual Schlesinger}}
\phantom{\begin{smallmatrix} \! \\ \! \\ \! \\ \! \\ \!\\ \!  \end{smallmatrix}}$ & $G^*$
\\ \hline
 $\sum \frac{A_i}{z-a_i}dz$ & 
$\phantom{\begin{smallmatrix} 1 \\ 1 \\ 1 \\ 1  \end{smallmatrix}}
\underset{\text{(classical Gaudin)}}{\text{Garnier}}
\phantom{\begin{smallmatrix} 1 \\ 1 \\ 1 \\ 1  \end{smallmatrix}}$
& Schlesinger & $G^n/G$
\\ \hline
  $\sum_1^3 \frac{A_i}{z-a_i}dz,\  \sl_2$ & & 
Painlev\'e VI & 
Fricke-Klein-Vogt
\\ \hline  
$\phantom{\begin{smallmatrix} . \\ . \\ . \\ . \\ . \end{smallmatrix}}
\underset
{\sl_2,\  C\text{ generic}}
{(A+Bz+Cz^2)dz}
\phantom{\begin{smallmatrix} . \\ . \\ . \\ . \\ . \end{smallmatrix}}$ 
& & 
Painlev\'e II & 
Flaschka--Newell
\\
  \vdots & \vdots & \vdots & \vdots
  \end{tabular}
\end{center}

The first row includes the Mi\v s\v cenko--Fomenko ``shift of argument'' integrable systems \cite{ms-shift} (for a recent review see 
\cite{bolsinov, FFT-shift}). 
For $B$ regular semsimple the corresponding isomonodromy systems fit into the class  considered by Jimbo--Miwa--M\^{o}ri--Sato \cite{JMMS}.
The corresponding (framed) wild character varieties are the standard dual Poisson Lie groups $G^*$ (the Poisson varieties underlying
the  Drinfeld--Jimbo quantum groups)--see \cite{smapg}.
Thus even simple wild character varieties 
seem important.
The underlying unframed/\hk spaces are  $\cL\spq_{\!q}\, T$ for 
symplectic leaves $\cL\subset G^*$. 
In this case the wild mapping class group coincides with the so-called quantum Weyl group \cite{bafi}.

The next row is the tame case on $\IP^1$: the isomonodromy system
is due to Schlesinger \cite{schles-icm1908}
and this led to the corresponding integrable system \cite{garn1919}.
The character varieties are the tame genus zero ones discussed above.
The next row is the case of this with $\sl_2$ and four poles, related to Painlev\'e VI and the Fricke--Klein--Vogt surfaces.

It turns out the first two rows are isomorphic, for general linear groups:
For the integrable systems see \cite{AHH-dual}, for the isomonodromy systems this is Harnad's duality \cite{Harn94}, for the character varieties this
follows from work on the Fourier--Laplace transform \cite{BJL81, malg-book} (summarised in \cite{cmqv} \S2), and for the full \hk metric one can upgrade this to a Nahm transform \cite{szabo-nahm, szabo-hk}. Their full structures as ``nonabelian Hodge spaces'' really are isomorphic.
For example we can take the Fricke--Klein--Vogt surface and 
ask how it arises explicitly 
on the irregular side in the first row,
in this case in terms of rank {\em three} bundles.
This can be done---see \cite{k2p} for the formulae.

More generally we expect that any tame character variety
of complex dimension two on the four-punctured sphere should be a Fricke--Klein--Vogt surface.
To test this take $G$ to be the simple group of type $G_2$.
It has dimension $14$ and has a special conjugacy class $\cC\subset G$ of dimension $6$ (a complex analogue of the $6$-sphere).
Then if we take $\cC_\infty\subset G$ to be a generic class the character variety 
\beq\label{eq: g2cv}
(\cC\fus{}\cC\fus{}\cC\fus{} \cC_\infty)\spq G\eeq
has dimension $3\cdot6+12-2\cdot 14 = 2$, and our expectation holds:
\begin{thm}[\cite{fricke}]
The tame $G_2(\IC)$ character variety \eqref {eq: g2cv} is isomorphic to a ``symmetric'' Fricke--Klein--Vogt surface \eqref{eq: fkv}, with $a=b=c$.
\end{thm}

\subsection{H3 surfaces}
More generally we expect all the cases of complex dimension two
are as follows (from the theory of Painlev\'e equations, and its extension by Sakai \cite{Sakai-CMP01}
related to connections on curves and thus \hk manifolds in \cite{quad}).
They are complete \hk four-manifolds, so are noncompact analogues of the K3 surfaces.
We call them H3 surfaces in honour of Higgs, Hitchin and Hodge.
The minimal (or ``standard'') representations
of these nonabelian Hodge spaces are as follows:
\begin{center}
 \begin{tabular}{| c | c | c | c | c | c |c | c | c | c | c| c| }
\hline
$\phantom{\bmx 1 \\ 1\emx}${Space} 
& $\wt E_8$ &  $\wt E_7 $ &  $\wt E_6$  &   $\wt D_4$ &  $\wt A_3=\wt D_3$ &  $\wt D_2$  & $\wt D_1$  & $ \wt D_0$ & $\wt A_2$ & $ \wt A_1$ & $\wt A_0$ \\
\hline
$\phantom{\bsmx 1 \\ 1 \\ 1\esmx}${Rank} 
&   $6$ &  $4$ &  $3$   &   $2$ &  $2$ &  $2$  & $2$  & $2$ & $2$ & $2$ & $2$ \\
\hline
$\phantom{\bmxnp 1 \\ 1 \emxnp}$
{Pole orders} 
$\phantom{\bmxnp 1 \\ 1 \emxnp}$
&   $111$ &  $111$ &  $111$   &   $1111$ &  $211$ &  $22$ & $2\wt2$ & $\wt2\,\wt2$ & $31$ & $4$ & $\wt 4$ \\
\hline
  \end{tabular}
\end{center}

All these are over the Riemann sphere;
{\em Rank} means the rank of the vector bundles
and  {\em Pole orders} means the orders of  the poles
of the connections/Higgs fields.
The tildes (\,$\wt\,$\,) on the pole order indicate a twisted irregular type 
(which here means having a nilpotent leading coefficient).
The standard representations of the first four cases are thus tame, and the other spaces have no known tame representations.
The cases $\wt E_8,\wt E_7$ require the residues to be in special adjoint orbits, determined by the affine Dynkin diagram 
(see \cite{quad}). More details/references are in \cite{ihptalk} \S3.2.

In each case the label indicates that the open part $\cM^*$ is diffeomorphic to the 
corresponding ALE or ALF \hk four-manifold, from
\cite{Kron.ale, hawking77, cherkishitchin, atiyah-hitchin, dancer93, houghton}. 
E.g. $\cM^*$ is the Eguchi--Hanson space for Painlev\'e II/$\wt A_1$ (as above), or $\IC^2$ for $\wt A_0$, or the Atiyah--Hitchin manifold for $\wt D_0$.
For the most part this matches the affine Weyl symmetry group of the Painlev\'e equation found by Okamoto \cite{Okamoto-dynkin}, but this viewpoint leads to a better understanding, since we can now generalise to higher dimensions.
(It is not clear how to generalise the ``perpendicular'' labelling preferred in  \cite{Okamoto-dynkin, Sakai-CMP01}.)

\subsection{Quiver varieties}

The higher dimensional analogue of 
Kronheimer's construction \cite{Kron.ale} 
of the 
ALE spaces are the additive/Nakajima quiver varieties 
\cite{nakaj-duke94}, 
and we can ask if some of them occur 
amongst the additive spaces $\cM^*$, generalising the special affine Dynkin graphs occuring in the H3 story, and the stars in 
\cite{CB-additiveDS}.
This is true and yields a theory of ``Dynkin diagrams'' for some of the nonabelian Hodge spaces.
In brief the quiver varieties arise (symplectically) as follows.
Given any graph $\Ga$ with nodes $I$, and an $I$-graded vector space $V$, there is a vector space $\Rep(\Ga,V)$ of
representations of $\Ga$ on $V$
(we view $\Ga$ as a doubled quiver---each edge denotes two oppositely oriented quiver edges).
This is a Hamiltonian $H$-space where $H=\Prod_I\GL(V_i)$.
Performing the symplectic reduction (at a central value $\la$ of the moment map) yields the quiver variety $\Rep(\Ga,V)\spq_{\!\la}\, H$.
 Now suppose $G=\GL_n$ and we have an irregular curve 
 $\bSi=(\IP^1,\infty,Q)$.
Let $\IC^n=\bigoplus_I V_i$ 
be the eigenspaces of $Q$, with eigenvalue 
$q_i\in z\IC[z]$ on $V_i$, so $H=\prod \GL(V_i)$ 
is the centraliser of $Q$ and $\cM^*\cong \cO_B\spq_{\!\La}\, H$ as above.
Let $\Ga(Q)$ be the ``fission graph'', 
with nodes $I$ and $\deg(q_i-q_j)-1$ edges between $i,j\in I$ 
(it was defined in an equivalent way, involving splaying/fission,
 in \cite{rsode} Apx. C).

 \begin{thm}[\cite{rsode, hi-ya-nslcase}]
 The orbit $\cO_B$ is isomorphic to $\Rep(\Ga(Q),V)$ as a Hamiltonian $H$-space, and consequently  $\cM^*$ is a Nakajima quiver variety.
 \end{thm}

The reduction by $H$ at $\La$ corresponds to gluing a leg (type $A$ Dynkin graph) onto each node of $\Ga(Q)$ to obtain a larger graph
$\wh \Ga(Q)$---we call such graphs ``supernova graphs'' as they generalise the stars
(cf. \cite{slims} Defn 9.1).
More generally  one can add some simple poles and still obtain that $\cM^*$ is a quiver variety (as in \cite{rsode}), giving more ``modular'' interpretations of certain  Nakajima quiver varieties (as moduli spaces of connections), making contact with  
\cite{CB-additiveDS} in the tame case.
E.g. from this  %
we can define spaces that should be the Hilbert schemes of the H3 surfaces (\cite{rsode} p.12, \cite{slims} \S11.4).

In the simply-laced case the fission graphs are exactly the complete $k$-partite graphs, determined by integer partitions with $k$ parts (see \cite{rsode}).
Each such graph can be ``read'' in terms of connections in $k+1$ different ways, giving isomorphic moduli spaces (both the open parts \cite{rsode, slims} and the full spaces 
\cite{cmqv}). 
\begin{figure}[h]
	\centering
	\input{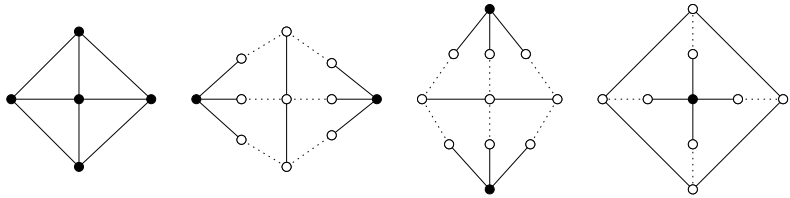}
	\caption{Four readings of 
$\Gamma(221)=2\cdot\Gamma(21)=1\cdot\Gamma(22)$, from \cite{rsode}.}
\label{fig: 221readings}
\end{figure}

This leads to the Kac--Moody Weyl group of the supernova graph (the ``global Weyl group'' of the corresponding irregular curve), and 
\cite{slims}
shows how it acts to give  
isomorphisms between the 
isomonodromy systems (generalising Harnad's duality), and between the integrable systems.
In turn considering the corresponding wild character varieties yields a new theory of multiplicative quiver varieties \cite{cmqv, even-euler}---we really can attach the whole nonabelian Hodge space to the graph.

\vspace{.2cm}

\noindent
{\Small
{\bf Acknowledgments.}
\noindent
The link to the $\wt D_0$-$\wt D_2$ ALF spaces came from 
discussions with Cherkis (in 2008 for $\wt D_2$, and
in the question session after the author's 2012 Banff talk
for $\wt D_0,\wt D_1$, in response to a question of Neitzke).
The notion of ``good'' came out of discussions with 
Yamakawa (cf. \cite{twcv}). 
}

\renewcommand{\baselinestretch}{1}              %
\normalsize
\bibliographystyle{amsplain}    \label{biby}
\bibliography{syr}

\vspace{0.5cm} 
\noindent
Math\'ematiques,
Universit\'e Paris-Sud, 91405 Orsay, France, Europe\\
Philip.Boalch@math.u-psud.fr

\end{document}